\documentclass{amsproc}

\newtheorem{theorem}{Theorem}[section]
\newtheorem{lemma}[theorem]{Lemma}
\newtheorem{proposition}[theorem]{Proposition}
\newtheorem{corollary}[theorem]{Corollary}

\theoremstyle{definition}
\newtheorem{definition}[theorem]{Definition}
\newtheorem{example}[theorem]{Example}

\theoremstyle{remark}

\numberwithin{equation}{section}

\newcommand{\moduli}{M}
\newcommand{\proj}{{\bf P}}
\newcommand{\rational}{{\bf Q}}
\newcommand{\real}{{\bf R}}
\newcommand{\ratmap}{--\!\!\!>}

\begin{document}

\title{Stable Pairs and Log Flips}

% author information
\author{Aaron Bertram}
\address{Department of Mathematics, University of Utah, Salt Lake City,
Utah 84112}
\email{bertram@math.utah.edu}
\thanks{Supported in part by a Sloan research fellowship}

\subjclass{14H10}
\date{December 12, 1996}

\maketitle

\section {Introduction.} 
Let $C$ be a Riemann surface of genus $g$, let $K_C$ be a canonical divisor
on $C$, and 
let $D$ be an effective divisor on $C$ of degree $d$. Two
basic rational maps arise in this context, namely the 
{\bf linear series map}:
$$\phi_{|K_C+D|}:C \ratmap |K_C+D|^*
\cong \proj^{d+g-2}$$
and the {\bf Serre correspondence}:
$$\psi_{|K_C+D|}:|K_C+D|^* \ratmap \moduli_C(2,D)$$
where $\moduli_C(2,D)$ is the moduli space for semistable bundles $E$ on $C$
satisfying rk$(E) = 2$ and det$(E) \cong {\mathcal O}_C(D)$. The latter 
arises when we consider hyperplanes in $|K_C+D|$ (via Serre duality) as 
lines in Ext$^1({\mathcal O}_C(D),{\mathcal O}_C)$ hence as extensions:
$$\epsilon:\ 0 \rightarrow {\mathcal O}_C \rightarrow E \rightarrow {\mathcal O}_C(D) \rightarrow 0
\ \ (\mbox{modulo equivalence and scalars}).$$
The following is an ancient result:
\begin{theorem}\label{old1}
If $d > 2$, then $\ \phi_{|K_C+D|}$ is an embedding. 
\end{theorem}
In case $d > 2$ and $d+g-2 = \ \mbox{dim}(|K_C+D|) > 1$, we will let
$$X := \ \mbox{bl}(|K_C+D|^*,C),$$
the blow-up of $|K_C+D|^*$ along the image of $C$ (which we will also denote by $C$).

In contrast, the Serre correspondence is a morphism 
only when $d \le 2$, as it is undefined at points of $C$
(and certain secant varieties) for 
larger values of $d$. The structure of the 
Serre correspondence for large $d$ (particularly  $d > 2g-2$,
in which case it is dominant) is, however, well understood now, due largely
to recent work of the author and Michael Thaddeus. The first goal 
of this paper (\S 2) is to explain Thaddeus' interpretation of the 
Serre correspondence as a sequence of simple birational maps of moduli 
spaces of stable pairs followed by a contraction.
We will discuss the birational maps at length, but let's start with an 
example of a contraction. This is another ancient result:

\begin{theorem}\label{old2}
If $d > 3$, then $C\subset |K_C+D|^*$ is cut out scheme-theoretically
by quadric hypersurfaces.  
\end{theorem}

If we let $H,E \subset X$ be the pull-back of a hyperplane divisor and 
an exceptional divisor,
respectively, then Theorem~\ref{old2} can be reinterpreted as the statement that
the linear series  $|2H - E|$ on $X$ is base-point free. As with
any base-point free
linear series, one obtains a well-defined contraction:
$$\gamma:X\rightarrow Y$$
(i.e. $\gamma$ satisfies 
$\gamma_*{\mathcal O}_X = {\mathcal O}_Y$).

\medskip 

For example, if  $g = 2$ and $d = 4$ then
$\gamma :X \rightarrow Y$ 
is generically a conic bundle over $Y \cong \moduli_C(2,{\mathcal O}_C(D)) \cong \proj^3$, 
and the blow-up followed
by the contraction resolve the Serre correspondence:
$$\begin{array}{ccccc}& & X
\\ & \swarrow & & \searrow \\ \psi_{|K_C+D|}: 
\proj^4 & & \ratmap & & \proj^3 \end{array}$$

In \S 2, we will show (following Thaddeus) how the moduli spaces of stable 
pairs produce   
varieties $Y_1 = |K_C+D|^*,X_1 = X,Y_2 = Y,X_2,Y_3,...,X_{[\frac{d-1}2]}$ 
and a diagram of contractions beginning with the blow-down 
and $\gamma$ and ending with a morphism to $\moduli_C(2,D)$ (the last morphism
is a contraction iff $d > 2g-2$): 

$$\begin{array}{cccccccccccccc}X & & & & X_2 & & & & & & & &
X_{[\frac{d-1}2]} \\
\downarrow & \stackrel \gamma \searrow & & \swarrow & & 
\searrow & & \swarrow & \cdots & 
\searrow & & \swarrow &  \downarrow\\
|K_C+D|^* & & Y & & & & Y_3 & & & & Y_{[\frac{d-1}2]} & & \moduli_C(2,D)
\end{array}$$

Just as the curve $C$ sits in $|K_C+D|^*$ via the linear series map,
there are natural inclusions of the $k$-th symmetric 
products $C_k$ of $C$ in each $Y_k$. (For $k > 1$, the 
$Y_k$ are singular along $C_k$). And just as the 
blow-down $X \rightarrow |K_C+D|^*$ contracts a projective bundle 
to $C$ and is an isomorphism elsewhere, so too is each of the 
contractions to a $Y_k$ an isomorphism off a projective 
bundle, which is contracted to $C_k$. 
Thus the diagram above decomposes the Serre correspondence as 
a sequence of ``irreducible'' birational maps $X_k \ratmap X_{k+1}$
followed by a morphism to $\moduli _C(2,D)$.

\bigskip

In the second part of the paper, we bring in the log minimal model 
program. The idea is to explain the diagram above
as a sequence of log flips starting from $X$.  
This will not provide a new {\bf proof} of the results
of \S 2, since the existence of log flips 
is only known in dimension $3$. Nevertheless,
the idea is to find an explanation for the existence of the diagram
without recourse to moduli spaces.
(I will use  vector bundle techniques, however, to prove the 
required properties of certain divisors on $X$).
The justification is the observation that from this 
point of view, it is not unreasonable to expect
more log flips than just those involving
moduli of stable pairs in case $d < 2g-2$. As motivation, recall the analogues
of Theorems~\ref{old1} and ~\ref{old2} for canonical curves
which are not explained by Thaddeus' stable pair construction. 

\begin{theorem}[Noether] If $C$ is not hyperelliptic, then
\ $\phi_{|K_C|}:C \ratmap |K_C|^*$ is an embedding.
\end{theorem}

\begin{theorem}[Petri] If $C$ is neither hyperelliptic nor 
trigonal nor a smooth plane quintic, then the image curve $C \subset |K_C|^*$ is 
cut out scheme-theoretically by quadric hypersurfaces. 
\end{theorem}

Section 2 begins with a general discussion of slope functions and 
stability.
The categories of triples (following Bradlow and Garcia-Prada \cite{BG}) and 
Bradlow pairs (following Thaddeus \cite{T1}) are investigated, 
together with 
their natural families of slope functions and stability conditions.
(A triple consists of two bundles 
plus a morphism. Triples are introduced here because the appearance of 
a parameter in the slope function has a very natural explanation in this 
context. Also, Bradlow pairs, which consist of a bundle plus a section,
can be viewed as a special case of triples.
See also \cite{B2}, \cite{RV} and in greatest generality, \cite{KN}, 
where moduli 
are constructed for 
pairs consisting of a vector bundle and a subspace of its space of 
sections --``Brill-Noether pairs''. A variant of this, where the subspace
comes along with a basis, is 
also considered in \cite{BDW}.)
Thaddeus' theorem on the existence of moduli spaces of stable  pairs 
in rank two is given, as well as the properties of these 
moduli spaces which are relevant to our discussion. 

Section 3 starts with some of the definitions and expectations
of the log minimal model program. Suitable
${\rational}$-divisors $F$ on $X$ are produced determinantally, and are 
proven to be log 
canonical, via a sequence of blow-ups of $X$ along smooth centers
first considered by the author in \cite{B1}.
Using these ${\rational}$-divisors, we are able to prove that all 
the simple birational maps $X_k \ratmap X_{k+1}$ produced by the diagram above 
are log flips, in the sense of the log minimal
model program. In contrast to 
\S 2,  several of the results of this section have not previously 
appeared elsewhere, and
are due to the author.

\medskip

{\bf Connections and Problems:} Three ``schools'' of algebraic 
geometry are represented in this paper -- curves in projective space,
moduli spaces and the minimal model program -- and my hope in writing it
is to stimulate more interaction among them. Concretely, I think there
are some interesting areas where a combination of ideas from the 
various schools may yield real progress. For example:

\medskip

{\bf Higher Dimension:} Stable pairs on surfaces have been studied 
(see \cite{HL} and LePotier's ``coherent systems'' \cite{L}), and 
the general problem of studying the behavior of GIT quotients 
(such as these) as the values of some parameters pass across a ``critical'' wall 
has been studied by \cite{DH} and \cite{T2}. However, to the best of 
my knowledge, no  
family of moduli spaces has been constructed which 
addresses the following systematically:

\medskip

{\bf Question:} Does there exist an analogous diagram of 
contractions and simple birational transformations for 
``sufficiently positive'' embeddings of higher dimensional varieties?
(And how positive is sufficiently positive?)

\medskip

{\bf Vanishing:} The existence of ``special'' log canonical divisors on 
$X$ would have some practical applications, even in the absence of 
a proof for the log minimal model program. There 
are general versions of Kodaira vanishing (see \cite{K}) which would 
apply in the presence of such divisors  to yield 
vanishings for the higher cohomologies of sheaves on $|K_C+D|^*$
of the form $I_C^k(n)$ (see \cite{B4}). 
From this point of view,
the following would yield concrete applications: 

\medskip

{\bf Question:} Given $d < 2g-2$, what is the minimum
value for a ratio $\frac ab$ of positive rational numbers
such that there exists an effective 
divisor $F \equiv aH - bE$ on $X$? 
And how large may $a$ be taken given $\frac ab$
so that $F$ is log canonical? 

\medskip

{\bf (Variant):} Same question if $C$ and/or $D$ are general in moduli.

\section{Stable Bradlow Pairs.} 
Let ${\mathcal A}$ be a category with a zero object in which kernels and cokernels
exist, as well as direct sums.
Let $S \subset Ob({\mathcal A})$ be a subset which is closed under direct sums. 

\begin{definition}
A function $\mu: S \rightarrow {\real}$
is called a {\it slope function} if 
for all short exact sequences $0 \rightarrow B \rightarrow A \rightarrow
C \rightarrow 0$ of elements of $S$, 
$$\mu(B) < \mu(A) \Leftrightarrow \mu(A) < \mu(C)\ \mbox{and}\  
\mu(B) = \mu(A) \Leftrightarrow \mu(A) = \mu(C)$$
\end{definition}

Given a slope function $\mu:S \rightarrow {\real}$:

\begin{definition} 
An object
$A\in S$ is called {\it stable} if $\mu(B) < \mu(A)$
whenever $B\in S$ and there exists an injection $B \hookrightarrow A$ 
other than the identity.
$A$ is called semistable if $\mu(B) \le \mu(A)$ above, 
strictly semistable if it is semistable but not stable, and
unstable if it is not semistable.
\end{definition}

\begin{definition} 
(a) If $A\in S$ is unstable, then a filtration:
$$0 \hookrightarrow A_1 \hookrightarrow A_2 \hookrightarrow ...
\hookrightarrow A_n = A$$
by elements of $S$ is called a Harder-Narasimhan filtration if 
the $A_i/A_{i-1}$ are all in $S$ and semistable and
$\mu(A_1) > \mu(A_2/A_1) > ... > \mu(A_n/A_{n-1})$.

\medskip

(b) If $A\in S$ is semistable, then a filtration
$$0 \hookrightarrow A_1 \hookrightarrow A_2 \hookrightarrow ...
\hookrightarrow A_n = A$$
is called a Jordan-H\"older filtration if the 
$A_i/A_{i-1}$ are all in $S$ and stable.
Given a Jordan-H\"older filtration of $A$, the object
$gr(A) := \oplus_{i=1}^nA_i/A_{i-1} \in S$ is called the associated graded of 
the filtration. Two Jordan-H\"older filtrations are called 
{\it s-equivalent} if 
their associated graded objects are isomorphic.
 
\end{definition}

\begin{example}{(Vector Bundles on $C$):}

$\bullet$ ${\mathcal A}$ is
the category of 
isomorphism classes of vector bundles on $C$. 

$\bullet$ $S = Ob({\mathcal A}) - \{0\}$.

$\bullet$ $\mu:S \rightarrow {\rational}$ is the usual 
slope function $\mu(E) = \ \mbox{deg}(E)/\mbox{rk}(E)$. Then:

(a) Harder-Narasimhan filtrations always exist (and up to 
isomorphism only depend upon 
the isomorphism class of $E$).

(b) Jordan-H\"older filtrations always exist, producing 
associated gradeds which only depend upon the isomorphism
class of $E$. In particular, $s$-equivalence
becomes an equivalence relation on isomorphism classes of 
semistable bundles.

(c) For fixed invariants $r$ (the rank) and either $d$ (the 
degree) or ${\mathcal O}_C(D)$ (the isomorphism class of the determinant),  
there are projective coarse moduli spaces
$\moduli_C(r,d)$ (respectively, $\moduli_C(r,D)$) for the functors ``families of semistable
vector bundles
modulo s-equivalence with the given invariants''. (See \cite{S} for details.)

\end{example}

The next example is due to Bradlow and Garcia-Prada
(\cite{BG}). 

\begin{example}{(Triples on $C$):} 

$\bullet$ ${\mathcal A}$ is the category 
of isomorphism classes of triples $(E,f,F)$, where $E,F$ are vector
bundles on $C$ and $f:E \rightarrow F$ is a homomorphism. A triple is called
{\it nontrivial} if $F\ne 0$ and {\it nondegenerate} if 
$f$ has maximal rank at some point. A morphism in this category is a pair
$(\alpha,\beta)$ consisting of morphisms 
$\alpha:E\rightarrow E'$ and $\beta:F \rightarrow F'$ in the category of 
vector bundles,
such that the following diagram commutes: 

$$\begin{array}{ccc}E & \stackrel\alpha\rightarrow  & E' \\
f\downarrow & & f'\downarrow \\ F & \stackrel \beta\rightarrow & F' \end{array}$$
Both $\alpha$ and $\beta$ need to be injective to make
$(\alpha,\beta): (E,f,F) \rightarrow (E',f',F')$
injective as a morphism of triples. Direct sums obviously
exist as the ``free'' sum of triples:
$(E,f,F) \oplus (E',f',F') = (E\oplus E',f\oplus f',F\oplus F')$.

$\bullet$ $S$ is the set of nontrivial triples.

$\bullet$ For each $\sigma \in {\real}$ and $(E,f,F) \in S$, let:
$$\mu_\sigma(E,f,F) = \frac{\mbox{deg}(E) + \ \mbox{deg}(F) + 
\sigma (\mbox{rk}(E) + \ \mbox{rk}(F))}{\mbox{rk}(F)}.$$

Note the asymmetry in the slope function!
One says $(E,f,F)$ is $\sigma$-stable if it is stable with respect to 
the slope function $\mu_\sigma$.
\end{example}

\begin{theorem}{(\cite{BG} Theorem 6.1)} For fixed invariants:
$$r_1 = \ \mbox{rk}(E), r_2 = \ \mbox{rk}(F), d_1 = \ \mbox{deg}(E),
d_2 = \ \mbox{deg}(F) \ \mbox{and}\ \sigma \in {\real}$$
a coarse moduli space exists for the functor
``families of non-degenerate $\sigma$-stable triples with the given invariants'',
which is moreover projective if $r_1 + r_2$ is relatively prime
to $d_1 + d_2$ and $\sigma$ is ``generic'' (see \cite{BG}). 
\end{theorem}

\noindent {\bf Explanation of the Parameter:} The idea is to relate 
stable triples $(E,f,F)$ on $C$ to (equivariantly) stable equivariant 
bundles $G$ on $C\times \proj^1$. (The action is the
automorphism group of $\proj^1$ acting on the second factor.) 
This is a consequence of K\"unneth, which gives an isomorphism 
$$\mbox{Hom}_{{\mathcal O}_C}(E,F) \cong \ 
\mbox{Ext}^1_{{\mathcal O}_{C\times \proj^1}}
(p^*E\otimes q^*{\mathcal O}_{\proj^1}(2),p^*F),$$
telling us to look for $G$ in the corresponding extension:
$$0 \rightarrow  p^*F \rightarrow G \rightarrow p^*E \otimes q^*{\mathcal O}_{\proj^1}(2)
\rightarrow 0.$$
(this technique is often called ``dimensional reduction''.)

The main point is now that ample line bundles on $C\times \proj^1$
are of the form $p^*L\otimes q^*M$,
and (an equivariant version of) Gieseker stability 
for bundles on $C\times \proj^1$ depends upon a parameter, 
namely the {\bf ratio} deg$(L)$/deg$(M)$. Stability of $G$ with respect to 
a given ratio translates into  $\sigma$-stability for 
triples for a fixed value of $\sigma$. 

\bigskip

We are most concerned with the following,
first considered by Bradlow in \cite{Br}.

\begin{example}{(Pairs on $C$):} 

$\bullet$ Restrict the category of 
triples to the objects: $({\mathcal O}_C,f,E)$ and 
$(0,0,E)$.  These objects we will call pairs, following the literature.
This full subcategory is closed under kernels and cokernels. It is not 
closed under arbitrary direct sums, but if a set of 
pairs is given, with the property that at most one of 
them is of the form  
$({\mathcal O}_C,f,E)$, then their direct sum does lie in the 
subcategory.
(This will be enough to construct associated gradeds for 
Jordan-H\"older filtrations!) 

$\bullet$ $S$ is the set of nontrivial pairs (i.e. $E \ne 0$) and a pair
of the form $({\mathcal O}_C,f,E)$ is nondegenerate if
and only if $f \ne 0$.

$\bullet$ The slope functions $\mu_\sigma$ are the same as for triples. 

\end{example}

With respect to this slope function, observe that 
a pair $({\mathcal O}_C,f,E)$ is 
$\sigma$-stable if and only if:

\medskip

(i)\  $\mu(F) < \mu(E) - \sigma(\frac 1{\mbox{\small rk}(F)} - 
\frac 1{\mbox{\small rk}(E)})$ 
for each $({\mathcal O}_C,g,F) \hookrightarrow ({\mathcal O}_C,f,E)$, and

\medskip

(ii)\ $\mu(F) < \mu(E) + \sigma(\frac 1{\mbox{\small rk}(E)})$ 
for each $(0,0,F) \hookrightarrow ({\mathcal O}_C,f,E)$.

\begin{theorem}{(\cite{T1} (1.1)-(1.19))} Fix invariants:
$$\mbox{rk}(E) = 2\ \mbox{and}\  
\mbox{det}(E) \cong {\mathcal O}_C(D),\ \mbox{with} \ 
\mbox{deg}(D) = d.\ \mbox{Then:}$$

(a) Harder-Narasimhan and Jordan-H\"older filtrations exist if $\sigma > 0$ and
yield a well-defined $s$-equivalence for $\sigma$-semi-stable pairs.

(b) For each $\sigma > 0$, a  projective variety $\moduli_C(2,D,\sigma)$ 
(abbreviated $\moduli_\sigma$)
coarsely represents the functor: 
``families of nondegenerate $\sigma$-semistable 
pairs $({\mathcal O}_C,f,E)$ modulo $s$-equivalence''. 
There is a universal family over the open locus parametrizing 
stable pairs, which is smooth and irreducible.

\end{theorem}
 
Full proofs of the properties listed below can be 
found in \cite{T1}. 

(Please note that our $\sigma$ differs from the $\sigma$
in \cite{T1} by a factor of $2$.)

\bigskip

\noindent {\bf Properties of Stable Pairs:} 
{\it Fix rk$(E) = 2$ and det$(E) \cong {\mathcal O}_C(D)$. 
Also assume that $g \ge 2$ (but see the note at the end of this section). Then:

\medskip

(a) There are no $\sigma$-semi-stable pairs if 
$\sigma < 0$ or $\sigma > d$.

\medskip

(b) There are always $0$-semi-stable pairs, though no $0$-stable pairs. 

\medskip

(c) $\moduli_{d}$ is a point.}

\bigskip

{\bf Proof:} If $\sigma < 0$, then by (ii) above, 
$(0,0,E) \hookrightarrow ({\mathcal O}_C,f,E)$
destabilizes any pair. If $\sigma > d$ and $({\mathcal O}_C,f,E)$
is given, let $L$ be the line-bundle image of ${\mathcal O}_C$ in $E$ 
with induced map $s:{\mathcal O}_C \rightarrow L$.
Then using (i) above, the pair is destabilized by the natural inclusion $({\mathcal O}_C,s,L) \hookrightarrow   
({\mathcal O}_C,f,E)$.

\medskip

If $\sigma = 0$, then conditions (i) and (ii) coincide, telling us that
$({\mathcal O}_C,f,E)$ is semistable if 
and only if $E$ is semistable, and that no pair is $0$-stable. Note that 
there is no Jordan-H\"older filtration.  

\medskip

If $\sigma = d$, then by the analysis in the proof of (a),
$({\mathcal O}_C,f,E)$ is $\sigma$-semi-stable 
if and only if $f:{\mathcal O}_C \rightarrow E$ has no zeroes,
and all such $({\mathcal O}_C,f,E)$ are $s$-equivalent, with associated graded 
$({\mathcal O}_C,\mbox{id},{\mathcal O}_C) \oplus  (0,0,{\mathcal O}_C(D))$. So 
there are no $\sigma$-stable pairs, and the 
moduli space is a point. 

\medskip

In contrast to the 
boundary cases presented here, the stable locus in $\moduli_\sigma$ will be 
nonempty if $0 < \sigma < d$.

\bigskip

\noindent {\bf Critical Points and Local Triviality:} 

\begin{center}

{\it Let $\Gamma = \{0 < c < d \ | \ 
c \equiv d \ (\mbox{mod}\ 2)\}$.}

\end{center}

{\it \noindent For each $\sigma$,  let $Z_\sigma \subset \moduli_\sigma$ be the locus of
$\sigma$-strictly-semistable pairs.

\medskip

(d) If $\sigma \not \in \Gamma$, then $Z_\sigma = \emptyset$ (i.e.
$\sigma$-semistable $\Rightarrow$
$\sigma$-stable).

\medskip

(e) If $c = d - 2n \in \Gamma$, then $Z_c \cong C_n$, the $n$-th symmetric product of $C$.

\medskip

(f) If $ I \subset (0,d) - \Gamma$ is an interval and $\sigma,\sigma' \in I$,
then $\moduli_\sigma \cong \moduli_{\sigma'}$.

\medskip

(g) Suppose that $c\in \Gamma$ and $c^- < c < c^+$ are real numbers
in the neighboring intervals of $(0,d) - \Gamma$. Then there are surjective 
morphisms:
$$\begin{array}{ccccc}\moduli_{c^-} &&&& \moduli_{c^+} \\
& \stackrel{f^-}\searrow & & \stackrel{f^+}\swarrow \\ && \moduli_c \end{array}$$
\noindent {\bf Key Point:} $f^-$ and $f^+$ are isomorphisms away from 
$Z_c \subset \moduli_c$ and projective bundles over $Z_c$. (The 
 projective bundles are identified in the proof).} 

\bigskip 

{\bf Proof:} If $\sigma > 0$ and if $F$ is the bundle in a destabilizing
subpair of $({\mathcal O}_C,f,E)$, then it is easy to see that $F$ is a line 
bundle. But if $\sigma \not \in \Gamma$, then the right side of 
(i) and (ii) are not integers, whereas $\mu(F)$ is an integer. So 
we cannot have equality. This proves (d).

\medskip

Suppose $c = d - 2n \in \Gamma$. Then a $c$-strictly semistable
pair $({\mathcal O}_C,f,E)$ has a subpair which is either isomorphic
to $({\mathcal O}_C,s,{\mathcal O}_C(A))$ (where $s$ is the tautological 
section, deg$(A) = n$, and ${\mathcal O}_C(A)$ is the image of ${\mathcal O}_C$ in $E$) or
else it is isomorphic to $(0,0,L)$, where deg$(L) = d - n$. But 
{\bf either} possibility forces the associated graded
for the pair $({\mathcal O}_C,f,E)$ to be
of the form $({\mathcal O}_C,s,{\mathcal O}_C(A)) \oplus (0,0,{\mathcal O}_C(D-A))$,
and these are parametrized by $C_n$.

\medskip

The stability conditions do not
change when $\sigma$ moves within an interval $I \subset (0,d) - \Gamma$
(again because $\mu(F) \in {\bf Z}$) so the moduli spaces are isomorphic
by the universal property of a coarse moduli space. 

\medskip

If $c\in \Gamma$, then apart from $Z_c$, the stability conditions do not change when
$c$ is replaced by $c^-$ or $c^+$, so the first part of the key point follows
as in the previous paragraph.

\medskip

Let $c = d - 2n$, and consider 
$({\mathcal O}_C,s,{\mathcal O}_C(A)) \oplus (0,0,{\mathcal O}_C(D-A)) \in Z_c$. Then
it follows that among all 
pairs with this associated graded, 
exactly those pairs of the form:
$$\begin{array}{ccccccccc}&&{\mathcal O}_C \\ && \downarrow & \searrow\\
0 & \rightarrow & {\mathcal O}_C(A) & \rightarrow & E & \rightarrow 
& {\mathcal O}_C(D-A) \rightarrow 0\end{array}$$
are $c^-$-stable, and these are parametrized by $|K_C + D - 2A|$, which has
dimension $d-2n+g-2$ (independent of $A$) since $d-2n > 0$.

\medskip

On the other hand, among all pairs with this associated graded,
exactly those pairs of the form:
$$\begin{array}{ccccccccc}&&&&&&{\mathcal O}_C \\ &&&&&\swarrow & \downarrow \\
0 & \rightarrow & {\mathcal O}_C(D-A) & \rightarrow & E & \rightarrow 
& {\mathcal O}_C(A) &\rightarrow & 0\end{array}$$
are $c^+$-stable. Such pairs are parametrized by  $\proj(V)$, where
$V$ sits in the long exact sequence:
$$H^0(C,{\mathcal O}_C(D-A)) \rightarrow V \rightarrow H^1(C,O_C(D-2A)) \rightarrow H^1(C,{\mathcal O}_C(D-A)).$$ 
(in fact, $\proj(V)$ is naturally isomorphic to $\proj(H^0(C,{\mathcal O}_C(D-A)
\otimes {\mathcal O}_A)^*))$.

\medskip

Thus the dimension of $\proj(V)$ is $n-1$, independent of $A$.
 
\bigskip

So there are a {\bf finite} number of moduli spaces $\moduli_\sigma$, linked by 
morphisms as in the following diagram:

$$\begin{array}{cccccccccccc} &X_2&&&& X_1 & & & & X_0\\ 
\cdots & & \searrow & & \swarrow & & \searrow & & \swarrow & & \searrow\\ 
& & & \moduli_{d-4} & & & & \moduli_{d-2} & & & & \moduli_d \\
\end{array}$$
where each $X_n \cong \moduli_{(d-2n)^-} \cong \moduli_{(d-2n-2)^+}$.

\medskip

Theorems ~\ref{old1} and ~\ref{old2} are embedded in this diagram because of:

\medskip

\noindent {\bf Large Values of $\sigma$:} 

\medskip

{\it (h) If $d > 0$, then $X_0 \cong |K_C+D|^*$.

\medskip

\hskip .22in  If $d > 2$, then $\moduli_{d-2} \cong |K_C+D|^*$ and $X_1 \cong X = 
\ \mbox{bl}(|K_C+D|^*,C)$. 

\medskip

\hskip .22in Moreover, the morphism $f^-: X_1 \rightarrow \moduli_{d-2}$ is the blow-down.

\medskip

\hskip .22in If $d > 4$, then $f^+:X_1 \rightarrow \moduli_{d-4}$
is the contraction $\gamma:X \rightarrow Y$.} 

\medskip

{\bf Proof:} A special case of the proof of (g) shows that $X_0 \cong |K_C+D|^*$.

\medskip

Another special case of the proof of (g) shows that 
$f^+:X_0 \rightarrow \moduli_{d-2}$ is an isomorphism, because the ``exceptional''
part of the map is a $\proj^0$-bundle(!)

\medskip

When (g) is applied to the map $f^-:X_1 \rightarrow \moduli_{d-2} = |K_C + D|^*$,
one discovers that the exceptional locus is a divisor, which is 
a projective bundle over $C$, hence $f^-$ is the blow-down.

\medskip

Finally, when (g) is applied to the map $f^+:X_1 \rightarrow \moduli_{d-4}$,
the exceptional set consists of lines spanned by two points of $C$
(i.e. the secant lines) which are contracted to points. This means that
the linear series which realizes $f^+$ must be a multiple of 
$|2H - E|$, so the fact that $f^+$ has connected fibers 
implies $f^+$ is equal to $\gamma$.

\bigskip

So we've got Theorem ~\ref{old1} and (a very precise) Theorem ~\ref{old2} when $d > 4$.
To see what happens for $d = 4$ from this point of view (for example, 
in the case $g = 2$ and $d = 4$ considered in the introduction), we need
to analyze:

\medskip

\noindent {\bf All Values of $\sigma$:} 

\medskip

{\it (i) If $n > 1$, then $X_n$ (if defined) is isomorphic to $X_1$ off codimension $2$.

\medskip

\hskip .2in The maps $f^+:X_n \rightarrow \moduli_{d-2n-2}$ are  multiples of
$|(n+1)H - nE|$.

\medskip

\hskip .2in The maps $f^-:X_n \rightarrow \moduli_{d-2n}$ are  
multiples of $|nH - (n-1)E|$.

\medskip

(j) If $\sigma$ is in the first interval of $(0,d) - \Gamma$, then 
there is a morphism:
$$f:\moduli_\sigma = X_{[\frac {d-1}2]} \rightarrow \moduli_C(2,D)$$
which is the contraction determined by high multiples of
$|dH - (d-2)E|$.}

\medskip

{\bf Proof:} The first part of (i) is a dimension count using
(g), which allows us to transfer linear series from $X_1 = X$ over
to each $X_n$. The reader is referred to \cite{T1}, where the 
ample cone is constructed for each $X_n$, the boundary of which gives 
properties (i) and (j). Notice in particular, that each $f^+$ and $f^-$
is a very simple contraction by property (g), but that the final 
map $f$ in (j) can have rather more complicated behavior, as in the 
example of the introduction.

\bigskip

Thus the mirror image of the diagram following property (g) gives us:
$$\begin{array}{cccccccccccc}X = X_1&&&&&&X_2 && & & X_{[\frac {d-1}2]}\\
\downarrow &&\searrow &&& \swarrow && \searrow &\cdots & \swarrow & \downarrow \\
|K_C+D|^*&&&&Y = \moduli_{d-4}&&&&& &\moduli_C(2,D)
\end{array}$$
which is the advertised generalization of Theorems ~\ref{old1} and ~\ref{old2}.

\bigskip

\noindent  Note: Most of this analysis also applies to curves of genus 
$0$ and $1$. 

\medskip

{\bf genus 1:} All properties (a)-(j) apply. The only difference between 
this and the general case occurs when $d$ is even, in which case
$\moduli_C(2,D)$ is isomorphic to $\proj^1$, rather than a point, as one
would expect by a dimension count. For example, if $d = 4$, then
$\gamma: X = \mbox{bl}(\proj^3,C) \rightarrow Y = \proj^1$ is the 
contraction determined by the pencil of quadrics vanishing along $C$.

\medskip

{\bf genus 0:} Properties (a)-(j) apply if $d$ is even and
$\moduli_C(2,D)$ is a point, corresponding to the vector bundle
${\mathcal O}_{\proj^1}(\frac d2) \oplus {\mathcal O}_{\proj^1}(\frac d2)$.
On the other hand, if $d = 2n + 1$, then $\moduli_\sigma = \emptyset$ if 
$\sigma < 1$, because all bundles are unstable. Other than 
this, which forces obvious changes to properties (b),(c),(g) and (j),
everything is as in the general case. Notice that in this case,
$\moduli_1$ is isomorphic to $\proj^n = (\proj^1)_n$, by property (e).
For example, when $d = 5$, then $\gamma: X = \mbox{bl}(\proj^3,C)
\rightarrow Y = \moduli_1 = \proj^2$ is the contraction determined by 
the web of quadrics vanishing along the twisted cubic. This contraction
is a $\proj^1$ bundle, a special case of the key point of property (g).

\section{Log Flips.} 
The goal of this section is to interpret the 
birational maps:
$$X = X_1 \ratmap X_2 \ratmap ... \ratmap X_{[\frac {d-1}2]}$$
as flips, in the sense of the minimal model program.
In fact, they are not flips, but rather their {\bf inverses} are flips
(at least initially),
in the traditional sense. While this is an interesting observation,
it is not the one I want to pursue, because the inverses point in the 
wrong direction, from the point of view of 
Theorems ~\ref{old1} and ~\ref{old2}. For example, with this interpretation, the first contraction
$\gamma:X\rightarrow Y$ is not a flipping contraction, but rather
contracts curves whose intersection with $K_X$ are positive. Such 
contractions are hard to understand, in general. Fortunately, the 
theory of log minimal models provides a means for turning the flips 
around, provided we can find suitable divisors on the $X_k$.
We begin with a quick tour of the parts of the log minimal model program
relevant to our discussion.  

\bigskip

Let $X$ be a smooth projective variety.

\begin{definition} A ${\rational}$-divisor on $X$ is a finite
sum of distinct prime divisors with rational coefficients. 
It is effective if all the 
coefficients are non-negative. Intersections with curves, self intersections
and numerical equivalence are all defined as with ordinary divisors.

Let $F$ be an effective ${\rational}$-divisor on $X$. If $F  = \sum \alpha_iF_i$,
then the support of $F$, denoted Supp$(F)$, is the union of the prime 
divisors $F_i$ which appear in $F$ with positive coefficients.

\end{definition}

\begin{definition} If $F$ is an effective divisor 
on $X$, then a {\bf log resolution} of $(X,F)$ is a morphism
$f:\widetilde X \rightarrow X$ with the property that $\widetilde X$ is 
smooth, and 
$\sum f\mbox{-exceptional divisors}\  + f_*^{-1}(\mbox{Supp}(F))$ is a 
normal crossings divisor.

\end{definition}

\noindent Note: If $f:Y \rightarrow X$ is any birational morphism
of smooth varieties, and 
if $D$ is a ${\rational}$-divisor on $X$,  
then $E_f$ will denote the sum of the $f$-exceptional divisors, and $f^*(D)$ and
$f_*^{-1}(D)$ will denote the total transform and the strict transform
of $D$ on $Y$, respectively. (They are well-defined by linearity.)

\medskip

Let $F = \sum \alpha_i F_i$ be an effective ${\rational}$-divisor.

\medskip

\begin{definition} $F$ is {\bf log canonical} if 
each coefficient $\alpha_i \le 1$, and there is a log resolution of 
$(X,F)$ with the property that all coefficients of the components of $E_f$
are at least $-1$ in the 
${\rational}$-divisor:
$$(K_{\widetilde X} -f^*K_X) +  (f_*^{-1}(F) - f^*(F))$$
\end{definition}

\noindent Note: This property is independent of the log resolution.

\begin{definition} Suppose that $B \subset X$ is a curve
(which we also identify with its image in $\mbox{H}_2(X,{\bf R})$). Then 
$B$ spans an {\bf extremal ray} of the cone of 
effective curves on $X$ if there is an element 
$\lambda \in \ \mbox{H}^2(X,{\bf R})$
such that:

(i) $\lambda(B) = 0$ and 

(ii) if $\beta \in \ \mbox{H}_2(X,{\bf R})$ is 
a limit of sums $\sum c_iB_i$ of curves with positive (real) coefficients,
then $\lambda(\beta) \ge 0$ with equality if and only if $\beta$ is 
a multiple of $B$.   
\end{definition}

\begin{definition} Suppose that $B$ spans an extremal ray and
$f:X \rightarrow Y$ is a morphism satisfying $f_*({\mathcal O}_X) = {\mathcal O}_Y$. 
If $B$ is
contained in a fiber of $f$, and if moreover every curve
contained in every fiber of $f$ is homologous to a (rational) multiple
of $B$, then $f$ (which is uniquely determined if it exists) is called
the {\bf extremal contraction} associated to $B$.
\end{definition}

Suppose that $F$ is a log-canonical divisor on $X$. A basic result
of the log minimal model program is the following (see \cite{CKM} and
\cite{Ketal}): 

\bigskip

\noindent {\bf Contraction Theorem:} {\it If $B\subset X$ 
spans an extremal ray 
and $B.(K_X + F) < 0$, then there is an extremal contraction 
$\gamma:X \rightarrow Y$ associated to $B$.}

\bigskip

A central question of the minimal model program is:

\bigskip

\noindent {\bf Do Log Flips Exist?:} Suppose the contraction $f:
X \rightarrow Y$ of the theorem is an isomorphism off 
codimension $2$ in $X$. Then 
does there exist a morphism $f^+:X^+ \rightarrow Y$ with the 
following properties:

\medskip

(a) $f^+$ is  
an isomorphism off codimension 2 in $X^+$. Let $(K_X + F)^+$ be the 
strict transform of $K_X+F$ in $X^+$.

\medskip

(b) If $B^+ \subset X^+$ is a curve lying in a fiber of $f^+$,
then $B^+.(K_X+F)^+ > 0$.

\medskip

(c) The singularities of $X^+$ (or rather, of  $(X^+,F^+)$) are
not too bad (for example, so that we can even define the intersections
$B^+.(K_X+F)^+$ in (b)).

\bigskip

When $F = \emptyset$ and the dimension of $X$ is $3$, then 
the affirmative answer to this question is a deep theorem of 
Mori (together with a definition of ``not too bad'', of course). 
The answer is also known to be yes for arbitrary $F$ and dimension 
$3$. The interested reader is urged to consult \cite{CKM}
and \cite{Ketal}, as well as Koll\'ar's notes in this Proceedings
for an introduction to the minimal 
and log minimal model programs and other applications. 
 
\bigskip

Next, we construct a 
morphism which will eventually be 
a log resolution. 

\bigskip

Let $M$ be a line bundle
on $C$, let $C_k$ be the $k$-th symmetric product of $C$,
and let $V = H^0(C,M)$. If $M$ has the following property:
$$(*)_k: \ \mbox{For all} \ D \in C_k, \ \mbox{dim}(H^0(C,M(-D)) = 
\mbox{dim}(V) - k$$
then each such divisor $D$ determines a $\proj^{k-1} \subset \proj(V)$,
which is called the span of $D$. Given  
that property $(*)_k$ holds, the $k$-th secant variety is:
$$\Sigma_k(C) = \bigcup _{D\in C_k} \ \mbox{span}(D) \subset \proj(V).$$
If $M$ satisfies $(*)_2$ (i.e. $M$ is very ample), 
let $X = \ \mbox{bl}(\proj(V),C)$ (as in \S 1).

\bigskip

\noindent Observation: ${\mathcal O}_C(K_C + D)$ satisfies $(*)_{d-1}$.
(Riemann-Roch!)

\medskip

The following construction 
blows up the secant varieties of $C$.

\begin{theorem}\label{logres} (\cite{B1} Theorem 1) 
(a) Suppose $n \ge 1$ and $M$ is a line bundle
with property $(*)_{2n}$. 
Then there is a birational morphism
$f:\widetilde X \rightarrow X$
which is a composition of the following blow-ups:

\medskip

$f^{(2)}: X^{(2)} \rightarrow X^{(1)} = X \ \mbox{blows up the strict transform
of $\Sigma_2(C)$}$,

\medskip

$f^{(3)}:X^{(3)} \rightarrow X^{(2)} \ \mbox{blows up the strict transform
of $\Sigma_3(C)$}$

\medskip

\hskip .5in $\vdots$

$f^{(n)}:\widetilde X = X^{(n)} \rightarrow X^{(n-1)}\ \mbox{blows up 
the strict transform
of $\Sigma_n(C)$}$

\medskip

Moreover, the strict transform of each $\Sigma_k(C)$ in $X^{(k-1)}$ is smooth 
and irreducible of dimension 
$2k-1$, transverse to all exceptional divisors, 
so in particular $\widetilde X$ is smooth. 

\medskip

For consistency, let $f^{(1)}: X \rightarrow \proj(V)$ be the blow-down.
Let $E^{(k)}$ be the strict transform in $\widetilde X$ of 
each $f^{(k)}$-exceptional divisor.
Then $E^{(1)} + ... + E^{(n)}$ is 
a normal crossings divisor on $\widetilde X$ with $n$ smooth components.

\medskip

If $M$ is a line bundle
that does not satisfy $(*)_2$,
let $\widetilde X = \proj(V)$. Then:

\medskip

(b) (Terracini recursiveness) 
Suppose $k \le n$ and 
$x\in \Sigma_k(C) - \Sigma_{k-1}(C)$. Then
the fiber 
$$(f^{(k)})^{-1}(x) \subset X^{(k)}$$
is naturally
isomorphic to $\proj(H^0(C,M(-2A)))$, where $A$ is the unique divisor of degree $k$
whose span contains $x$.
Moreover, the fiber 
$$f^{-1}(x) \subset E^{(k)} \subset \widetilde X$$
is isomorphic to $\widetilde X_A$, 
the variety obtained by applying (a) of the Theorem to
the line bundle $M(-2A)$. 

\medskip

(c) If $g \ge 2$ and if $M ={\mathcal O}_C(K_C+D)$, then there is a  
{\bf morphism} 
$$\widetilde \psi_{|K_C+D|}: \widetilde X \rightarrow \moduli_C(2,D)$$
which extends $\psi_{|K_C+D|}$.
When restricted to a fiber
$f^{-1}(x)$ of part (b), $\widetilde \psi_{|K_C+D|}$ agrees with 
$\widetilde \psi_{|K_C+D-2A|}$ (and this property determines 
$\widetilde \psi_{|K_C+D|}$ uniquely!) 

\end{theorem}

For the proof, see \cite{B1}. Notice that
parts (a) and (b) make 
no reference to moduli, hence generalize to, for example, canonical 
embeddings, where condition $(*)_d$ is equivalent to the nonexistence
of $g^1_d$'s. As for part (c), the idea is to construct a vector bundle
on $C \times \widetilde X$ by a sequence of elementary modifications
of the bundle (constructed from the universal 
extension) on $C \times |K_C+D|^*$
along the exceptional divisors for each $f^{(k)}$, and to use 
this bundle to get the map to moduli. 

\medskip

In fact, though, the proof really constructs families of nondegenerate pairs
$({\mathcal O}_C,f,E)$ parametrized by the $\widetilde X^{(k)}$ (in all genera)
with the following property. For every $y\in \widetilde X$ and every 
$\sigma \in [0,d)$
(or $[1,d)$ if $g = 0$ and $d$ is odd) there is an $X^{(k)}$ such that 
the image of $y$ in $X^{(k)}$ parametrizes a $\sigma$-semistable pair.
Thus, for each $\sigma$, there is a natural morphism:
$$\psi_\sigma:\widetilde X \rightarrow \moduli_\sigma.$$

\bigskip

Now we construct log-canonical divisors on $X = \ \mbox{bl}(|K_C+D|^*,C)$.

\medskip

\noindent {\bf Linear Algebra Construction:} Given any vector bundle $F$ on $C$, 
the cup product gives rise to a linear map:
$$c: \ \mbox{Ext}^1(F(D),F) \rightarrow 
\ \mbox{Hom}(H^0(C,F(D)),H^1(C,F))$$

Also, the summand ${\mathcal O}_C \hookrightarrow F \otimes F^*$ 
produces an inclusion of vector spaces:
$$\iota: \ \mbox{Ext}^1({\mathcal O}_C(D),{\mathcal O}_C) \hookrightarrow 
\ \mbox{Ext}^1(F(D),F)$$

One can think of the composition $c\circ \iota$
pointwise as follows. Given
$$\epsilon:\ 0 \rightarrow {\mathcal O}_C \rightarrow E \rightarrow {\mathcal O}_C(D) 
\rightarrow 0$$
one tensors each term by $F$, and $c(\iota(\epsilon))$ is the 
connecting homomorphism:
$$c(\iota(\epsilon)) = \delta: H^0(C,F(D)) \rightarrow H^1(C,F).$$

When we lift $c\circ \iota$ to a map of trivial bundles 
on $|K_C + D|^*$, 
it determines a matrix $M(F)$ of linear forms on $|K_C+D|^*$ via:
$${\mathcal O}_{|K_C+D|^*}(-1) \rightarrow \ \mbox{Hom}(H^0(C,F(D)),H^1(C,F))
\otimes {\mathcal O}_{|K_C+D|^*}.$$

\begin{proposition}\label{linalg} (a) For each $0 < k \le \frac d2$, there is a 
nonempty open subset $U \subset \ \mbox{Pic}^{k-(g-1)}(C)$ such that
$$(*)\ \ L \in U \Rightarrow \ h^0(C,L^{-1}(D)) = d-k \ \mbox{and}\ 
h^1(C,L^{-1}) = k.$$

If $L \in U$, choose a basis for $H^0(C,L^{-1}(D))$, and let 
$I = (i_1,...,i_k)$ be a mutiindex with $1 \le i_1 < ... < i_k \le d-k$.
Then $I$ determines a $k\times k$ minor $M_I(L^{-1})$ 
(choosing columns $i_1,...,i_k$ from the matrix $M(L^{-1})$) yielding a
divisor on $\widetilde X$:
$$D_{L,I} \in 
|kH - (k-1)E^{(1)} - (k-2)E^{(2)} - ... - E^{(k-1)}|.$$
(i.e. the generic multiplicity of $\mbox{det}(M_I(L^{-1}))$
along $\Sigma_i(C)$ is at least $k-i$).  

\medskip

Finally, if we let $V_k$ be the sub-linear-series spanned by the $D_{L,I}$,
then $V_k$ is base-point-free (and independent of choices of basis).

\medskip

(b) Suppose $g > 0$. Then for each $0 < l \le d$, there is a nonempty 
open subset $U \subset \moduli_C(2,l-(2g-2))$ such that 
$$(*)\ \ F \in U \Rightarrow \ h^0(C,F^{-1}(D)) = 2d - l\ \mbox{and}\ 
h^1(C,F^{-1}) = l.$$

If $F \in U$, choose $J = (j_1,...,j_l)$ such that
$1 \le j_1 < ... < j_l \le d-l$ and the minor $M_J(F^{-1})$ as in (a). Then
det$(M_J(F^{-1}))$ determines a divisor:
$$D_{F,J} \in |lH - (l-2)E^{(1)} - (l-4)E^{(2)} - ... |,$$
and the sub-linear-series $W_l$
spanned by the $D_{F,J}$ is base-point-free.
\end{proposition}

{\bf Proof:} The values for
$h^0$ and $h^1$ in $(*)$ are generic in Pic and $\moduli_C(2,*)$ 
respectively. $U$ is 
an intersection of
two nonempty open subsets.

For the next part, the following observation is crucial. Given
an effective  divisor $A$ on $C$, an extension
$\epsilon \in \ \mbox{Ext}^1({\mathcal O}_C(D),{\mathcal O}_C)$ determines
a point $\overline \epsilon \in \ \mbox{span}(A) \subset |K_C+D|^*$
if and only if the extension splits when pushed forward:
$$\begin{array}{ccccccccc}0 & \rightarrow & {\mathcal O}_C & \rightarrow 
& E & \rightarrow & {\mathcal O}_C(D) & \rightarrow & 0 \\
&&\downarrow & \swarrow  \\ &&{\mathcal O}_C(A) \end{array}$$

Now suppose $\overline \epsilon \in \Sigma_i(C)$, so is in the 
span of some divisor $A$ of degree $i < k$. It follows
(tensoring the inclusion ${\mathcal O}_C(D-A) \hookrightarrow E$ 
by $L^{-1}$) that 
$H^0(C,L^{-1}(D-A)) \subset \ \mbox{ker}(c(\iota(\epsilon)))$, which 
by Riemann-Roch has dimension at least $(d-k)-i$. Thus the rank of each 
$M_I(L^{-1})$ is at most $i$, and therefore its determinant has 
multiplicity at least $k-i$ at $\overline \epsilon$, from which the 
linear series computation in (a) follows. The linear series in (b) is 
computed similarly.
 
\medskip

We prove base-point-freeness first when $g=0$ and $d = 2n+1$. 
Given $k$, if  
$y  \in  \widetilde X - 
(E^{(1)} \cup ...\cup E^{(k-1)})$, then 
the bundle $E$ associated to $\overline \epsilon = y$ 
is isomorphic
to ${\mathcal O}_{\proj^1}(m) \oplus {\mathcal O}_{\proj^1}(d-m)$ where
$d-m > m \ge k$. This 
is because of the crucial observation. 
It follows that $h^0(C,E(-k-1)) = d-2k$, so some
$M_I({\mathcal O}_{\proj^1}(-k-1))$ has full rank at $\overline \epsilon$,
and thus $y$  is not a base point. If  
$y\in E^{(i)}$ for some (minimal) $i < k$, then Theorem ~\ref{logres} (b) allows us to 
place $y$ in a fiber over $\Sigma_i(C)$ isomorphic to  
$\widetilde X_A$ for some 
divisor $A$ satisfying deg$(A) = i$. The restriction
to this fiber of $V_k$ is identified with the linear series $V_{k-i}$
under the isomorphism with $\widetilde X_A$, and so we can conclude 
base-point-freeness by induction.

\medskip

The proof of base-point-freeness is similar in general. Suppose  
$d = 2n+1$ or $2n+2$, so ${\mathcal O}_C(K_C+D)$ satisfies $(*)_{2n}$. 
Given $k$ (or $l = 2k$ or $2k-1$), first consider the points
$y \in \widetilde X - (E^{(1)} \cup ... \cup E^{(k-1)})$. 
If $y \in E^{(k)} \cup ... \cup E^{(n)}$, then
 we can
find an $L \in U$ (or $F\in U$) such that 
$h^0(C,E\otimes L^{-1}) = d - 2k$ (or $h^0(C,E\otimes F^{-1}) = 2d - 2l$)
because the bundle $E$ associated to $\overline \epsilon = f(y)$
fits in an exact sequence:
$0 \rightarrow {\mathcal O}_C(D-A) \rightarrow E 
\rightarrow {\mathcal O}_C(A) \rightarrow 0$
($k \le \ \mbox{deg}(A) \le n$).
(This is a consequence of the crucial observation.)
If $y$ does not
lie in an exceptional divisor, then the bundle $E$ associated to the extension 
$\overline \epsilon = y$ is semistable
(Theorem ~\ref{logres} (c) or the crucial observation) in which  case 
the same fact about $h^0(C,E\otimes L^{-1})$ (or $h^0(C,E\otimes F^{-1})$)
is a standard result,
for example, see \cite {B3}, Lemma 3.6. Once this is 
achieved, 
one has base-point-freeness off the exceptional 
divisors $E^{(1)} \cup ... \cup E^{(k-1)}$ and 
the points of these exceptional divisors 
are treated by induction using Theorem ~\ref{logres} (b) and the same identification
of $V_k$ with $V_{k-i}$ (or $W_l$ with $W_{l-2i}$) as above.

\bigskip

There are two exceptional cases (given in detail below) where the linear series
($V_k$ in case (a) and $W_l$ in case (b)) are trivial, which occur when
$\Sigma_n(C)$ is a divisor and $k$ (or $l$) is maximal. In all other cases,
we can use Bertini to find smooth members of the linear series which meet
all the exceptional divisors $E^{(1)},...,E^{(n)}$ transversally.

\bigskip

\noindent {\bf Exceptional Cases:} (a) Suppose $g = 0$ and $d = 2n+2 = 2k$. Then
there is only one $k\times k$ matrix $M_I({\mathcal O}_{\proj^1}(-k-1))$,
and $V_k$ has only one element. (So since it is base-point-free, it has
to be trivial!) 
With a suitable choice of basis, 
$M_I({\mathcal O}_{\proj^1}(-k-1))$) is the 
standard square matrix:

\medskip

$$\left(\begin{array}{ccccc}z_1 & z_2 & z_3 & \cdots & z_k \\
z_2 & z_3 & z_4 & \cdots & z_{k+1} \\
\vdots & \vdots & \vdots & & \vdots \\
z_k & z_{k+1} & z_{k+2} & \cdots & z_{d-1}\end{array}\right)$$ 

\medskip

\noindent whose $2\times 2$ minors cut out the rational normal curve
$C \subset \proj^{d-2} = |K_C+D|^*$ (see Proposition 9.7 in \cite{H}).
Its determinant cuts out $\Sigma_n(C) \subset \proj^{d-2}$. 
The linear series $V_k$ is trivial because here
$E^{(n)} \equiv nH - (n-1)E^{(1)} - ... - 2E^{(n-1)}$ on $\widetilde X$.
Notice that although we cannot use Bertini to find smooth divisors
in $V_k$, we can use Theorem ~\ref{logres} (a) to conclude that $E^{(n)}$ itself
is smooth and meets the other exceptional divisors transversally.

\medskip

(b) Suppose $g = 1$, and $d = l = 2n+1$. Then there is 
one $l\times l$ matrix $M_J(F^{-1})$ for each stable bundle $F$ 
of rank $2$ and degree $d = d-(2g-2)$. There is no reason a priori 
why this determines a trivial linear series,
however, as in Exception (a), one computes that 
$\Sigma_n(C) \subset |K_C+D|^*$ is a divisor, of degree $d$,
which must therefore be the zero locus of each determinant, 
and $E^{(n)} \equiv dH - (d-2)E^{(1)} - ... - 3E^{(n-1)}$. 
(The degree can be computed using Lemma 2.5 (Chapter VIII)
from \cite{ACGH}). Again, we will use the fact that $E^{(n)}$ 
is smooth, intersecting the other
exceptional divisors transversally.

\bigskip

\noindent {\bf Remarks:} If $l$ is even, the conditions
$h^0(C,F^{-1}(D)) = 2d-l$ and $h^0(C,F^{-1}) = l$ may not be
independent of the choice of a representative $F$
for a semistable point in $\moduli_C(2,l-(2g-2))$. However, if these properties are 
true for the associated graded, then they hold for all representatives,
as is easily checked. 
Moreover, if we let 
$l = 2k$, then the split bundles determine an inclusion of 
linear series: $V_k \cdot V_k \subseteq W_l$. (In genus $0$ and $1$, this is
an equality!) 

\bigskip

Next, we use the linear series to find:

\bigskip

\noindent {\bf Some Log Canonical Divisors on $X$:} Let $F = f_*A$, where:

\medskip

\noindent {\bf Genus 0:} (a) If $d = 2n+1$, then $A \in W_{2n}$ 
is a general member.

\medskip

(b) If $d = 2n+2$, then $A = E^{(n)} + A'$, where 
$A' \in V_{n}$ is a general member.

\medskip

\noindent {\bf Genus 1:} (a) If $d = 2n+1$, then $A = E^{(n)}$.

\medskip

(b) If $d = 2n+2$, then $A \in W_d$ is a general member.

\medskip

\noindent {\bf Genus $\ge$ 2:} $A \in W_d$ is general.

\bigskip

Alternatively, one can think of $F$ in each case 
as the strict transform in $X$ of a
hypersurface in $|K_C+D|^*$ (highly singular along the secant varieties). 
However, when we think of $F$ as the push-forward of a divisor $A$
on $\widetilde X$, then
the following becomes almost immediate.

\bigskip

\noindent {\bf Claim:} In all the cases above, $f:\widetilde X \rightarrow X$ is 
a log resolution of $(X,F)$ and $F$ is a log canonical divisor.

\medskip

{\bf Proof:} By Theorem ~\ref{logres} (a), all the $f$-exceptional divisors
are smooth with normal crossings. In each case, the  strict transform of 
the support of $F$ 
is the support of $A$, which is a sum of smooth divisors which 
intersect all others with normal 
crossings, either by Bertini or Theorem ~\ref{logres}(a) again. So $f$ is a 
log resolution of $(X,F)$.   

\medskip

Since each blow-up $f^{(k)}$ was along a smooth center trasnverse to all 
exceptional divisors, it follows that the coefficient of $E^{(k)}$ in
$K_{\widetilde X} - f^*K_X$ is the codimension of 
$\Sigma_k(C)$ in $|K_C+D|^*$ minus $1$, a  consequence
of Riemann-Hurwitz. It also follows that since we constructed
$F$ as $f_*A$, the coefficient
of $E^{(k)}$ in $f_*^{-1}F - f^*F = A-f^*F$ is the (negative of the)
generic multiplicity
of $F$ along the strict transform of $\Sigma_k(C)$ in $X$, which 
is computed directly from the linear series in which $A$ 
lies.
This is the information we need to check that $F$ is log canonical.
The computations in genus $0$ are left to the reader. Here 
is the data for genus $\ge 1$:

\medskip

$\bullet$ \ codimension of $\Sigma_k(C)$ in $|K_C+D|^*$:\ \ $d + g - 2k -1$.

\medskip

$\bullet$ \ multiplicity of $F$ along $E^{(k)}$:\ \ $d - 2k$.

\medskip

Since it follows that the coefficient of each $E^{(k)}$ in 
$(K_{\widetilde X} - f^*K_X)
+ (A - f^*F)$ is  $g-2$, we see that $F$ is log canonical. 

\medskip

We can (and need to!) do a little better when $g \ge 2$ 
if we use ${\rational}$-divisors. If $p,q$ are positive integers, let
$(W_d)^{p}$ be the linear series spanned by products of $p$ elements of $W_d$,
and given a smooth element $G \in (W_d)^p$ (this linear series is 
base-point free), consider $F' = \frac 1qf_*(G)$. This is 
not only numerically equivalent to $\frac pq F$ (as is easy to see), but
all coefficients of the $E^{(k)}$ in the expression 
$f_*^{-1} F' - f^*F'$
are $\frac pq$ times the corresponding coefficients for $F$. 
We will abuse notation and say that this divisor is a member of $\frac pq F$, keeping in 
mind that if $p > q$, then the literal ${\rational}$-divisor $\frac pq F$ 
cannot be log canonical, by definition, while a member 
constructed in this way might be log canonical. In fact, 
if $d > 4$, then 
$$\left(\frac{d+g-5}{d-4}\right)F \ \mbox{has a log canonical member}$$   
by the data above (keep in mind that $E^{(1)}$ is not
$f$-exceptional).
 
\bigskip

Now we will relate these log canonical divisors to 
the diagram at the end of \S 2 constructed by stable pairs. Namely, recall that 
whenever $d > 2k$, there was a diagram:
$$\begin{array}{ccccc} X_{k-1} & & & & X_k \\
&\stackrel{f^+}\searrow \ & & \stackrel{f^-}\swarrow \\ & & 
\moduli_{d-2k}\ \ \ \ \end{array}$$

Moreover, $f^-$ and $f^+$ are obviously extremal ray contractions
since each contracts a projective bundle over $C_k$ and 
dim$(H_2(X,\real)) = 2$. (Take any
curve in a projective-space fiber to span the extremal ray.)
Finally, each contraction is an isomorphism off of codimension $2$.

\bigskip

\begin{proposition}\label{logflip} If $k = 2$ or $d > 2g -2$ and $k$ is arbitrary,
then the diagram above
is a log flip for (the strict transform on $X_{k-1}$ of) 
$K_X + \left(\frac {d+g-5}{d-4}\right)F$.
\end{proposition}

{\bf Proof:} We need to show: (a) the member of 
$\left(\frac {d+g-5}{d-4}\right)F$ 
constructed as above is log canonical
on each $X_{k}$ (this will certainly suffice for the condition
``not too bad'' in the ``definition'' of log flips), and
(b) If $B \subset X_{k-1}$ and 
$B^+ \subset X_k$ are curves spanning extremal rays, then 
$B.(K_X+\left(\frac {d+g-5}{d-4}\right)F) < 0$
and $B^+.(K_X+\left(\frac {d+g-5}{d-4}\right)F) > 0$.

\bigskip

We prove (b) first.
Recall (property (i) from \S 2) that the map $f^+$ is a multiple of 
$|kH - (k-1)E|$. Thus, if $B \subset X_{k-1}$ is an 
extremal ray, then $B.(H - \frac{k-1}kE) = 0$. Moreover, 
$|(k-1)H - (k-2)E|$ is nef on $X_{k-1}$, and $B$ is 
{\bf not} contracted in this linear series, so it follows that
$B.(H-\frac{k-2}{k-1}E) > 0$, and $B.E > 0$.

\medskip

\noindent From the data: 

\medskip

$K_X \equiv -(d+g-1)H + (d+g-4)E$, and 

\medskip

$F \equiv dH - (d-2)E$,

\medskip

\noindent we get $K_X+\left(\frac{d+g-5}{d-4}\right)F \equiv \frac {4g-4}{d-4}(H - \frac{d+2g-6}{4g-4}E)$,
from which it follows that its intersection with $B$ is negative
 when $k = 2$ or 
$d > 2g-2$. Moreover, the mirror 
image of this argument shows that if $B^+ \subset X_k$ is 
an extremal ray for $f^-$, then its intersection
with $K_X+\left(\frac{d+g-5}{d-4}\right)F$ is positive in the same cases.

\medskip

So the only thing left to see is the fact that 
the member of $\left(\frac {d+g-5}{d-4}\right)F$ we constructed is 
log canonical
on all $X_k$, not just $X = X_1$, as was shown earlier. In fact, I claim
a stronger result, which will explain all the maps as log flips:

\begin{lemma} If $d > 2k$, then a general member of
$\left(\frac{d+g-2k-1}{d-2k}\right)F$ is
log canonical on $X_{k-1}$.
\end{lemma}

{\bf Proof:} The construction is as before, pushing down  
a general element of $(W_d)^{d+g-2k-1}$ and dividing by
$d-2k$.
After $k-1$ elementary modifications, the 
proof of Theorem ~\ref{logres}(c) (see \cite{B1}) produces a family
of $\sigma = n-2k+1$-stable pairs on $C$ parametrized by $X^{(k-1)}$, hence
a morphism $\psi^{(k-1)}:X^{(k-1)} \rightarrow X_{k-1}$ since
$X_{k-1} = \moduli_\sigma$. Moreover,
the morphism $\psi_\sigma : \widetilde X \rightarrow X_{k-1}$ factors 
through $\psi^{(k-1)}$  via the composition of blow-downs
$f_{k-1}:\widetilde X \rightarrow X^{(k-1)}$ from Theorem ~\ref{logres}(a). 

\medskip

One checks that if
$G \in (W_d)^p$ (for any $p > 0$), then 
$(f_{k-1})_*G$ descends to a divisor on $X_{k-1}$. 
Thus when we log resolve $(X_{k-1},\frac pq F)$ by the map
$\psi_\sigma$, then only the 
exceptional divisors $E^{(k)}$ and above appear with a nonzero 
coefficient in $\psi_\sigma^*\frac pq F - (\psi_\sigma)_*^{-1}(\frac pq F)$,
and those, it is easy to see, appear with the same coefficients 
as in the  earlier computation. The lemma immediately follows. 

\medskip

\begin{corollary}\label{logcan} Each rational map $X_{k-1} \ratmap X_{k}$ is
a log flip.
\end{corollary}

{\bf Proof:} Using the log canonical divisor on $X_{k-1}$ from the Lemma:
$$K_X + \left(\frac{d+g-2k-1}{d-2k}\right)F \equiv
\frac {k(2g-2)}{d-2k}\left(H - \frac{(2k-1)(2g-4) + d-2}{k(2g-2)}E\right)$$
has negative intersection with $B$ and positive intersection
with $B^+$ (as in the proof of Proposition ~\ref{logflip}, keeping in 
mind the fact that $d > 2k$). 

\bigskip

\noindent {\bf Final Remark:} I have split off Proposition ~\ref{logflip}
from Corollary ~\ref{logcan} (which is in a sense more powerful!) to point
out a curious fact. Namely, if $d > 2g-2$, which is precisely  when
$\psi_{|K_C+D|}:|K_C+D|^* \ratmap 
\moduli_C(2,D)$ is dominant, then we can construct a single ${\rational}$-divisor on $X$ 
for which all the maps $X_{k-1} \ratmap X_k$ simultaneously become log flips.
When $d \le 2g-2$, however, one needs to tailor the divisor to 
the variety $X_{k-1}$ and in fact it seems that no single ${\rational}$-divisor 
on $X$ will be log canonical and have the desired intersection properties with 
all the extremal rays. (At least the linear algebra construction does not 
produce such a divisor.) 

\medskip

\noindent {\bf Acknowledgements:} I would like to thank Michael Thaddeus 
and the referee for their careful reading and useful comments on an 
earlier version of this paper.

University of Utah, Salt Lake City, UT 84112

{\it email address:} bertram@math.utah.edu

\end{document}